\documentclass[10pt]{article}
\usepackage{fullpage}

\usepackage{times}
\usepackage[latin1]{inputenc} 
\usepackage[T1]{fontenc} 
\usepackage{url}
\usepackage{boxedminipage}
\usepackage{amsmath,amstext,amsfonts,amssymb}
\usepackage{listings}
\usepackage{color}
 \newcommand{\crlibm}{\texttt{CRlibm}}

\newcommand{\R}{\mathbb{R}}
\definecolor{gray}{gray}{0.5}
\definecolor{dkgray}{gray}{0.3}

\lstdefinelanguage{gappa}{
  comment=[l]\#,
  otherkeywords={->,/\\, ?, $, \{, \}, float },      
  morekeywords={in},
  basicstyle=\small\ttfamily,
  keywordstyle=\small\ttfamily,
}

\title{Certifying floating-point implementations using Gappa}
\author{Florent de Dinechin,
        Christoph Lauter,
        and Guillaume Melquiond}

\begin{document}

\maketitle 

\begin{abstract}
  High confidence in floating-point programs requires
  proving numerical properties of final and intermediate values. One
  may need to guarantee that a value stays within some range, or that
  the error relative to some ideal value is well bounded. Such work
  may require several lines of proof for each line of code, and will
  usually be broken by the smallest change to the code (e.g. for
  maintenance or optimization purpose).  Certifying these programs by hand
  is therefore very tedious and error-prone.  This article discusses the use of
  the Gappa proof assistant in this context.  Gappa has two main
  advantages over previous approaches: Its input format is very close
  to the actual C code to validate, and it automates error evaluation
  and propagation using interval arithmetic.  Besides, it can be used
  to incrementally prove complex mathematical properties pertaining to
  the C code. Yet it does not require any specific knowledge about
  automatic theorem proving, and thus is accessible to a wide
  community. Moreover, Gappa may generate a formal proof of the
  results that can be checked independently by a lower-level proof
  assistant like Coq, hence providing an even higher confidence in the
  certification of the numerical code. The article demonstrates the
  use of this tool on a real-size example, an elementary function with
  correctly rounded output.

\end{abstract}


\section{Introduction}
\label{sect:intro}

Floating-point
(FP) arithmetic was designed to help developing software
handling real numbers.  However, FP numbers are only an
approximation to the real numbers. A novice programmer may incorrectly
assume that FP numbers possess all the basic properties of the real numbers,
for instance associativity of the addition, and waste time fighting
the associated subtle bugs. Having been bitten once, he may forever
stay wary of FP computing as something that cannot be trusted. As many
safety-critical systems rely on floating-point arithmetic, the question of the
confidence that one can have in such systems is of paramount
importance, all the more as floating-point hardware, long available in
mainstream processors, is now also increasingly designed into embedded
systems.

This question was addressed in 1985 by the IEEE-754
standard for floating-point arithmetic~\cite{IEEE754}. This standard
defines common floating-point formats (single and double precision),
but it also precisely specifies the behavior of the basic operators
$+$, $-$, $\times$, $\div$, and $\sqrt{~~}$.  In the rounding mode to the
nearest, these operators shall return the \emph{correctly rounded}
result, uniquely defined as the floating-point number closest to the
exact mathematical value (in case of a tie, the number returned is the
one with the even mantissa). The standard also defines three
\emph{directed rounding} modes (towards $+\infty$, towards $-\infty$,
and towards $0$) with similar correct rounding requirements on the
operators.

The adoption and widespread use of this standard have increased the
numerical quality and portability of floating-point code. It has
improved confidence in such code and allowed construction of proofs of
numerical behavior~\cite{Goldberg91}. Directed rounding modes are also
the key to enable efficient \emph{interval arithmetic}~\cite{Moore66},
a general technique to obtain validated numerical results.

This article is related to the IEEE-754 standard in two ways.
Firstly, it discusses the issue of proving properties of numerical
code, building upon the properties specified by this standard.
Secondly, one of our case studies will be the proof of the implementation of an
elementary function (the logarithm) returning correctly-rounded
results. Such correctly-rounded implementations will contribute to a better
floating-point standard, where the numerical quality of elementary
functions matches that of the basic operators.

Elementary functions were left out of the IEEE-754 standard in 1985 in
part because the correct rounding property is much more difficult to
guarantee than for the basic arithmetic operators. The problem will be
exposed in detail in the sequel. In short,
the correctness of the implementation of a correctly rounded function
requires the \emph{a priori} knowledge of a bound on the overall error
of evaluating $f(x)$. Moreover, this bound should be tight, as a loose
bound
negatively impacts performance \cite{DinErshGast2005,DinLauMul2007:log}. 

This article describes an approach to machine-checkable proofs of such
tight bounds, which is both interactive and easy to manage, yet much
safer than a hand-written proof. It applies to error bounds as well as
range bounds. Our approach is not restricted to the validation of
elementary functions. It currently applies to any straight-line
floating-point program of reasonable size (up to several hundreds of operations). 

The novelty here is the use of a tool that transforms a
high-level description of the proof into a machine-checkable version,
in contrast to previous work by Harrison~\cite{Har97,harrison-fmcad2000} who directly
described the proof of the implementation of some functions
in all the low-level details.  The Gappa approach is more concise and
more flexible in the case of a subsequent change to the code. More
importantly, it is accessible to people outside the formal proof
community.

This article is organized as follows. Next section describes in detail
the challenges posed by automatic computation of tight bounds on ranges and errors.
Section~\ref{sec:gappa} describes the Gappa tool. Section~\ref{sec:tutorial}
gives an overview on the techniques for proving an
elementary function using Gappa and give an extensive example of the interactive
construction of the proof.

\section{Proving properties of floating-point code}
\label{sec:bound}

\subsection{Floating-point numbers are not real numbers}
We have already mentioned that floating-point (FP) numbers do not possess basic
properties of real numbers. The following FP code sequence, due to Dekker \cite{Dekker71}, illustrates this:

 \begin{lstlisting}[language={C},firstnumber=1,label=lst:Fast2Sum,caption=The Fast2Sum algorithm.]
   s = a+b;  
   r = b-(s-a);           
\end{lstlisting}

This sequence consists only of three operations. The first one
computes the FP sum of the two numbers $\mathtt{a}$ and $\mathtt{b}$. The
second one would always return $\mathtt{b}$ and the third one $0$, if this FP sum was exact. As the numbers
here are FP numbers, however, the sum is often inexact, because of the rounding. In IEEE-754
arithmetic with round-to-nearest, under certain conditions, this
algorithm computes in $\mathtt{r}$ the error committed by this first
rounding.  In other words, it ensures that
$\mathtt{r}+\mathtt{s}=\mathtt{a}+\mathtt{b}$ in addition to the fact
that $\mathtt{s}$ is the FP number closest to $\mathtt{a}+\mathtt{b}$.
The Fast2Sum algorithm provides us with an exact representation of the
sum of two FP numbers as a pair of FP numbers, a very useful operation.

This example illustrates an important point, which pervades all of this
article: FP numbers may be an approximation of the reals that fails to
ensure basic properties of the reals, but they are also a very
well-defined set of rational numbers, which have other well-defined
properties, upon which it is possible to build mathematical proofs such as the proof
of the Fast2Sum algorithm.

Let us come back to the condition under which the Fast2Sum algorithm
works. This condition is that the exponent of \texttt{a} is larger
than or equal to that of \texttt{b}, which will be true for instance
if $|\mathtt{a}|\ge |\mathtt{b}|$. If one is to use this algorithm,
one has first to prove that this condition is true.  Note that
alternatives to the Fast2Sum exist for the case when one is unable to
prove this condition. The version by Knuth \cite{Knuth97} requires $6$
operations instead of $3$. Here, being able to prove the condition, which is a
property on values of the code, will result in better performance.

The proof of the properties of the Fast2Sum sequence (three FP
operations) requires several pages \cite{Dekker71}, and is indeed
currently out of reach of the Gappa tool, basically because it can not
be reduced to manipulating ranges and errors. This is not a problem,
since this algorithm has already been proven using machine checkers
\cite{DauRidThe01}. We consider it as a generic building-block of
larger floating-point programs, and the focus of our approach is to
automate the proof of such larger programs. In the case of the
Fast2Sum, this means proving the condition.

Let us now introduce the class of larger FP programs that we have
targeted in this work:  implementations of elementary functions

\subsection{Elementary functions}

Current floating-point implementations of elementary functions
\cite{Gal86,Tang91,StoTan99,Markstein2001,Muller2006} have several
features that make their proof challenging:

\begin{itemize}
\item The code size is too large to be safely proven by hand. In the
  first versions of the \crlibm{} project, the complete paper proof of a
  single function required tens of pages. It is difficult to trust
  such proofs.
\item The code is optimized for performance, making extensive use of
  floating-point tricks such as the Fast2Sum above. As a consequence,
  classical tools of real analysis cannot be
  straightforwardly  applied. Very often, considering the same operations on
  real numbers would be simply meaningless.
\item The code is bound to evolve for optimization purpose, because
  better algorithms may be found, but also because the processor
  technology evolves. Such changes will impose to rewrite the proof,
  which is both tedious and error-prone.
\item Much of the knowledge required to prove error bounds on the code
  is implicit or hidden, be it behind the semantics of the programming
  language (which defines implicit parenthesing, for examples), or in
  the various approximations made. Therefore, the mere translation of
  a piece of code into a set of mathematical variables that represent
  the values manipulated by this code is tedious and error-prone if
  done by hand.
\end{itemize}

Fortunately, FP elementary function implementations also
have pleasant features that make their proof tractable:

\begin{itemize}
\item There is a clear and simple definition of the mathematical object that the
  floating-point code is supposed to approximate. This will not always
  be the case of e.g. scientific simulation code. 
\item The code size is small enough to be tractable, typically less than a
  hundred floating-point operations.
\item The control flow is simple, consisting mostly of straight-line code with a few tests but no
  loops.
\end{itemize}

The following elaborates on these features.

\subsubsection{A primer on elementary function evaluation}

The evaluation of an elementary function is
classically \cite{Markstein2001,Muller2006} performed by a polynomial approximation valid on a small
interval only. A \emph{range reduction} step brings the input number
$x$ into this small interval, and a \emph{reconstruction} step builds
the final result out of the results of both previous steps. For
example, the logarithm may use as a range reduction the errorless
decomposition of $x$ into its mantissa $m$ and exponent $E$: $x = m
\cdot 2^E$. It may then evaluate the logarithm of the mantissa, and
the reconstruction consists in evaluating $\log(x) = \log(m) + E \cdot
\log(2)$.  Note that current implementations typically involve several
layered steps of range reduction and reconstruction. With current
processor technology, efficient implementations~\cite{Gal86, Tang91,
  Priest97} rely on large tables of precomputed values. See the books
by Muller~\cite{Muller2006} or Markstein~\cite{Markstein2000} for recent
surveys on the subject.

In the previous logarithm example, the range reduction was exact, but the
reconstruction involved a multiplication by the irrational $\log(2)$,
and was therefore necessarily approximate.  This is not always the case. For example,
for trigonometric functions, the range reduction involves subtracting
multiples of the irrational $\pi/2$, and will be inexact, whereas the
reconstruction step consists in changing the sign depending on the
quadrant, which is exact in floating-point arithmetic.

It should not come as a surprise that either range reduction or
reconstruction are inexact. Indeed, FP numbers are rational numbers,
but for most elementary functions, it can be proven that, with the
exception of a few values, the image of a rational is irrational.
Therefore, in an implementation, one is bound, at some point, to
manipulate numbers which are approximations of irrational numbers. 

This introduces another issue which is especially relevant to
elementary function implementation. One wants to obtain a
double-precision FP result which is a good approximation to the
mathematical result, the latter being an irrational most of the times.
For this purpose, one needs to evaluate an approximation of this
irrational to a precision better than that of the FP format. 

\subsubsection{Reaching better-than-double precision}

Better-than-double precision  is typically attained thanks to \emph{double-extended}
arithmetic on processors that support it, or \emph{double-double}
arithmetic, where a number is held as the unevaluated sum of two
doubles, just as the 8-digit decimal number $3.8541942\cdot 10^{1}$  may be represented by the
 unevaluated sum of two 4-digit numbers $3.854\cdot 10^{1}+1.942\cdot
 10^{-3}$.  
 Well-known and well-proven algorithms exist for manipulating
 double-double numbers~\cite{Dekker71,Knuth97},  the simplest of which is the Fast2Sum already introduced. 

However these algorithms are
 costly, as each operation on double-double numbers requires several FP
 operations. In this article, we will consider implementations based
 on double-double arithmetic, because they are more challenging, but
 Gappa handles double-extended arithmetic equally well.

\subsection{Floating-point errors, but not only}

The evaluation of any mathematical function entails
two main sources of errors.
\begin{itemize}
\item Approximation errors (also called methodical errors), such as
  the error of approximating a function with a polynomial. One may
  have a mathematical bound for them (given by a Taylor formula for
  instance), or one may have to compute such a bound using numerics
  \cite{ChevLaut2007QSIC}, for example if the polynomial has been
  computed using Remez algorithm.

\item Rounding errors, produced by most floating-point operations of
  the code. 
\end{itemize}

The distinction between both types of errors is sometimes arbitrary:
for example, the error due to rounding the polynomial coefficients to
floating point numbers is usually included in the approximation error of the
polynomial. The same holds for the rounding of table values, which is
accounted for more accurately  as approximation error than  as rounding error. This point is mentioned here because a lack of
accuracy in the definition of the various errors involved in a given
code may lead to one of them being forgotten.

\subsection{The point with efficient code}

Efficient code is especially difficult to analyze and prove because of all
the techniques and tricks used by expert programmers.

For instance, many floating-point operations are exact, and the
experienced developer of floating-point code will try to use them.
Examples include multiplication by a power of two, subtraction of
numbers of similar magnitude thanks to Sterbenz' Lemma~\cite{Ste74},
exact addition and exact multiplication algorithms (returning a
double-double), multiplication of a small integer by a floating-point
number whose mantissa ends with enough zeroes, etc. 

The expert programmer will also do his best to avoid computing more
accurately than strictly needed. He will remove from the computation
some operations that are expected not to improve the accuracy of the
result by much. This can be expressed as an additional approximation.
However, it soon becomes difficult to know what is an approximation to
what, especially as the computations are re-parenthesized to maximize
floating-point accuracy.

To illustrate the resulting code obfuscation, let us introduce the
piece of code that will serve as a running example along this
article.

\subsection{Example: A double-double polynomial evaluation}
\label{sec:troislignes} Listing \ref{lst:troislignes} is an extract of
the code of a sine function in the \crlibm{} library. These three
lines compute the value of an odd polynomial $p(y) = y + s_3\times y^3
+ s_5\times y^5 + s_7\times y^7$ close to the Taylor approximation of
the sine (its degree-1 coefficient is equal to $1$). In our algorithm,
the reduced argument $y$ is ideally obtained by subtracting to the FP
input $x$ an integer multiple of $\pi/256$. As a consequence $y \in
[-\pi/512, \pi/512] \ \subset \ [-2^{-7},
2^{-7}]$.

However, as $y$ is irrational, the implementation of this range
reduction has to return a number more accurate than a double, 
otherwise there is no hope of achieving an accuracy of the sine that matches
floating-point double precision. In our implementation, range
reduction therefore returns a double-double $\mathtt{yh}+\mathtt{yl}$.

To minimise the number of operations, Horner rule is used for the
polynomial evaluation: $p(y) = y + y^3\times(s_3 + y^2 \times (s_5 +
y^2 \times s_7))$.  For a double-double input
$y=\mathtt{yh}+\mathtt{yl}$, the expression to compute is thus
$(\mathtt{yh}+\mathtt{yl}) + (\mathtt{yh}+\mathtt{yl})^3\times(s_3 +
(\mathtt{yh}+\mathtt{yl})^2\times(s_5 +
(\mathtt{yh}+\mathtt{yl})^2\times s_7)).  $

The actual code uses an approximation to this expression: the
computation will be accurate enough if all the Horner steps
except the last one are computed in double-precision. Thus, $y_l$ will be
neglected for these iterations, and coefficients $s_3$ to $s_7$ will
be stored as double-precision numbers noted $\mathtt{s3}$, $\mathtt{s5}$
and $\mathtt{s7}$. The previous expression becomes:
\[(\mathtt{yh}+\mathtt{yl}) + 
\mathtt{yh}^3\times(\mathtt{s3} + \mathtt{yh}^2\times(\mathtt{s5} + \mathtt{yh}^2\times\mathtt{s7})). \]

However, if this expression is computed as parenthesized above, it
will not be very accurate. Specifically, the floating-point addition
$\mathtt{yh}+\mathtt{yl}$ will (by definition of a double-double)
return $\mathtt{yh}$, so the information held by $\mathtt{yl}$ will be
completely lost. Fortunately, the rest of the Horner evaluation also
has much smaller magnitude than $\mathtt{yh}$ (this is deduced from $y\in[-2^{-7}, 2^{-7}]$, therefore $y^3\in[-2^{-21}, 2^{-21}]$). 
The following parenthesing is therefore much more accurate:
\[\mathtt{yh} + \left(\mathtt{yl} + \mathtt{yh} \times
\mathtt{yh}^2\times(\mathtt{s3} + \mathtt{yh}^2\times(\mathtt{s5} + \mathtt{yh}^2\times\mathtt{s7}))\right). \]

In this last version of the expression, only the leftmost addition must be
accurate: we will use a Fast2Sum (which as we saw is an exact addition of two doubles
returning a double-double). The other operations use the native --- and
therefore fast --- double-precision arithmetic. We obtain the code of
Listing~\ref{lst:troislignes}.

 \begin{lstlisting}[language={C},firstnumber=1,label=lst:troislignes,caption=Three lines of C]
  yh2 = yh*yh;					   
  ts = yh2 * (s3 + yh2*(s5 + yh2*s7));	   
  Fast2Sum(sh,sl,   yh,  yl + yh*ts);	           
\end{lstlisting}

To sum up, this code implements the evaluation of a polynomial with
many layers of approximation. For instance, variable \texttt{yh2}
approximates $y^2$ through the following layers:
\begin{itemize}
\item $y$ was approximated by $\mathtt{yh}+\mathtt{yl}$ with the relative accuracy $\epsilon_{\mathrm{argred}}$
\item $\mathtt{yh}+\mathtt{yl}$  is approximated by $\mathtt{yh}$  in most of the computation,
\item $\mathtt{yh}^2$ is approximated by   \texttt{yh2=yh*yh}, with a floating-point rounding error. 
\end{itemize}

In addition, the polynomial is an approximation to the sine function,
with a relative error bound of $\epsilon_{\mathrm{approx}}$ which is
supposed known (how it was obtained it is out of the scope of this
paper \cite{ChevLaut2007QSIC}).

Thus, the difficulty of evaluating a tight bound on an elementary
function implementation is to combine all these errors without
forgetting any of them, and without using overly pessimistic bounds
when combining several sources of errors. The typical trade-off here
will be that a tight bound requires considerable more work than a
loose bound (and its proof might inspire considerably less
confidence). Some readers may get an idea of this trade-off by relating
each intermediate value with its error to confidence intervals, and
propagating these errors using interval arithmetic. In many cases, a
tighter error will be obtained by splitting confidence intervals into
several cases, and treating them separately, at the expense of an
explosion of the number of cases. This is one of the tasks that Gappa
will helpfully automate.

\subsection{Previous and related work}

We have not yet explained why a tight error bound is required in order
to obtain a correctly rounded implementation. This question is
surveyed in \cite{DinLauMul2007:log}. To sum it up, an error bound is
needed to guarantee correct rounding, and the tighter the bound, the
more efficient the implementation.  A related problem is that of
proving the behaviour of \emph{interval} elementary functions
\cite{DinMaid2006:RNC}.  In this case, a bound is required to ensure
that the interval returned by the function contains the image of the
input interval. A loose bound here means returning a larger interval
than possible, and hence useless interval bloat. In both case, the
tighter the bound, the better the implementation.

As a summary, proofs written for version of the \crlibm{} project up to
versions 0.8~\cite{crlibmweb} are typically composed of several pages of paper proof and
several pages of supporting Maple for a few lines of code. This
provides an excellent documentation and helps maintaining the code, but
experience has consistently shown that such proofs are extremely
error-prone. Implementing the error computation in Maple was a first
step towards the automation of this process, but if it helps avoiding
computation mistakes, it does not prevent methodological mistakes. 
Gappa was designed, among other objectives, in order to fill this void.

There has been other attempts of assisted proofs of elementary
functions or similar floating-point code. The pure formal proof
approach of Harrison \cite{Har97,harrison-fmcad2000,harrison2003-sqrt}
goes deeper than the Gappa approach, as it accounts for approximation
errors. However it is accessible only to experts of formal proofs, and
fragile in case of a change to the code. The approach of Krämer
\emph{et al} \cite{filib98,filib2005} relies on operator overloading
and does not provide a formal proof.



\section{The Gappa tool}
\label{sec:gappa}

Gappa\footnote{\url{http://lipforge.ens-lyon.fr/www/gappa/}}
extends the interval arithmetic para\-digm to the
field of numerical code certification~\cite{DauMel04,MelPio05}. Given the
description of a logical property involving the bounds of mathematical
expressions, the tool tries to prove the validity of this property.
When the property contains unbounded expressions, the tool computes
bounding ranges such that the property holds. For instance, the incomplete property ``$x + 1 \in [2,3] \Rightarrow x \in [?,?]$''
can be input to Gappa. The tool answers that $[1,2]$ is a range of
the expression $x$ such that the whole property holds. 

Once Gappa has reached the
stage where it considers the property to be valid, it generates a
formal proof that can be checked by an independent proof assistant.
This proof is
completely independent of Gappa and its validity does not depend on
Gappa's own validity. It can be mechanically verified by an external
proof checker and included in bigger formal developments.

\subsection{Floating-point considerations}

Section~\ref{sec:tutorial} will give examples of Gappa's syntax and
show that Gappa can be applied to mathematical expressions much more
complex than just $x+1$, and in particular to floating-point
approximation of elementary functions. This requires describing
floating-point arithmetic expressions within Gappa.

Gappa only manipulates expressions on real numbers. In the property $x +
1 \in [2,3]$, $x$ is just a universally-quantified real number and the
operator $+$ is the usual addition on real numbers $\R$. Floating-point
arithmetic is expressed through the use of ``rounding operators'':
functions from $\R$ to $\R$ that associate to a real number $x$ its
rounded value $\circ(x)$ in a specific format. These operators are
sufficient to express properties of code relying on most floating-point
or fixed-point arithmetics.

Verifying that a computed value is close to an ideal value can now be
done by computing an enclosure of the error between these two values.
For example, the property ``$x \in [1,2] \Rightarrow \circ(\circ(2 \times
x) - 1) - (2 \times x - 1) \in [?,?]$'' expresses the absolute
error caused during the floating-point computation of the following
numerical code:
\begin{lstlisting}[language=C,firstnumber=1]
float x = ...;
assert(1 <= x && x <= 2);
float y = 2 * x - 1;
\end{lstlisting}

Infinities and NaNs (Not-a-Numbers) are not an implicit part of this formalism:
The rounding operators
return a real value and there is no upper bound on the magnitude of
the floating-point numbers. This means that NaNs and
overflows will not be generated nor propagated as they would in
IEEE-754 arithmetic. However, one may use Gappa to prove very useful
properties, for instance that overflows, or NaNs due to some division
by $0$,  cannot occur in a given code: This can be expressed in terms
of intervals. What one cannot prove are properties depending on the
correct propagation of infinities and NaNs in the code.

\subsection{Proving properties using intervals}

Thanks to the inclusion property of interval arithmetic, if $x$ is an
element of $[0,3]$ and $y$ an element of $[1,2]$, then $x + y$ is an
element of the interval sum $[0,3] + [1,2] = [1,5]$. This technique
based on interval evaluation can be applied to any expression on real
numbers. That is how Gappa computes the enclosures requested by the user.

Interval arithmetic is not restricted to this role though. Indeed the
interval sum $[0,3] + [1,2] = [1,5]$ do not only give bounds on $x + y$,
it can also be seen as a proof of $x + y \in [1,5]$. Such a computation
can be formally included as an hypothesis of the theorem on the enclosure
of the sum of two real numbers. This method is known as computational
reflexivity~\cite{Bou97} and allows for the proofs to be machine-%
checkable. That is how the formal proofs generated by Gappa can be
checked independently without requiring any human interaction with a
proof assistant.

Such ``computable'' theorems are available for the Coq~\cite{HueKahPau04}
and HOL Light~\cite{Har00} proof assistants. Previous work~\cite{DauMelMun05}
on using interval arithmetic for proving numerical theorems has shown
that a similar approach can be applied for the PVS~\cite{OwrRusSha92} proof
assistant. As long as a proof checker is able to do basic computations
on integers, the theorems Gappa relies on could be provided. As a
consequence, the output of Gappa can be targeted to a wide range of
formal certification frameworks, if needed.

\subsection{Other computable predicates}

Enclosures are not the only useful predicates. As intervals are
connected subsets of the real numbers, they are not powerful enough to
prove some properties on discrete sets like floating-point numbers or
integers. So Gappa handles other classes of predicates for an expression
$x$:
\begin{eqnarray*}
\mathtt{FIX}(x, e) & \equiv & \exists m \in \mathbb{Z},\ x = m \cdot 2^e\\
\mathtt{FLT}(x, p) & \equiv & \exists m, e \in \mathbb{Z},\ x = m \cdot 2^e \land |m| < 2^p
\end{eqnarray*}

As with intervals, Gappa can compute with these new predicates. For
example, $\mathtt{FIX}(x,e_x) \land \mathtt(y,e_y) \Rightarrow
\mathtt{FIX}(x+y,\min(e_x,e_y))$. These predicates are especially useful
to detect real numbers exactly representable in a given format. In
particular, Gappa uses them to find rounded operations that can safely
be ignored because they do not contribute to the global rounding error.
Let us consider the floating-point subtraction of two floating-point
numbers $x \in [3.2,3.3]$ and $y \in [1.4,1.8]$. Note that Sterbenz'
Lemma is not sufficient to prove that the subtraction is actually exact,
as $\frac{3.3}{1.4} > 2$. Gappa is, however, able to automatically prove
that $\circ(x - y)$ is equal to $x - y$.

As $x$ and $y$ are floating-point numbers, Gappa first proves that they
can be represented with $24$ bits each (assuming single precision
arithmetic). As $x$ is bigger than $3.2$, it can then deduce it is a
multiple of $2^{-22}$, or $\mathtt{FIX}(x,-22)$. Similarly, it proves $\mathtt{FIX}(y,-23)$. The
property $\mathtt{FIX}(x-y,-23)$ comes then naturally. By computing with
intervals, Gappa also proves that $|x-y|$ is bounded by $1.9$. A
consequence of these last two properties is $\mathtt{FLT}(x-y,24)$: only
$24$ bits are needed to represent $x-y$. So $x-y$ is representable by a
single-precision floating-point number.

There are also some specialized predicates for enclosures. The following
one expresses the range of an expression $u$ with respect to another
expression $v$:
\begin{eqnarray*}
\mathtt{REL}(u,v,[a,b]) & \equiv & -1 < a \land \exists \epsilon \in [a,b],\ u = v \cdot (1 + \epsilon)
\end{eqnarray*}
This predicate is seemingly equivalent to the enclosure of the relative
error $\frac{u - v}{v}$. It allows, however, to simplify proofs, as the
error can now be manipulated even when $v$ is potentially zero. For
example, the relative rounding error of a floating-point addition
vanishes on subnormal numbers (including zero) and is bounded elsewhere,
so the following property holds: $\mathtt{REL}(\circ(x + y), x + y,
[-2^{-53},2^{-53}])$ when rounding to nearest in double precision.

\subsection{Gappa's engine} \label{sec:engine}

Because basic interval evaluations do not keep track of correlations
between expressions sharing the same terms, some computed ranges may be
too wide to be useful. This is especially true when bounding errors. For
example, when bounding the absolute error $\circ(a) - b$ between an
approximation $\circ(a)$ and an exact value $b$. The simplest option is to first compute the
ranges of $\circ(a)$ and $b$ separately and then subtract them. 
However, first rewriting the expression as $(\circ(a) - a) + (a - b)$
and then bounding $\circ(a) - a$ (a simple rounding error) and $a - b$
separately before adding their ranges usually gives a much tighter
result. These rules are inspired by techniques developers usually apply
by hand in order to certify their numerical applications.

Gappa includes a database of such rewriting rules (section
\ref{sec:hints} shows how the user may expand this database). The tool
applies them automatically, so that it can bound expressions along
various evaluation paths. Since any of these resulting intervals
encloses the initial expression, their intersection does, too. Gappa
keeps track of the paths that lead to the tightest interval
intersection and discards the others, so as to reduce the size of the
final proof. It may happen that the resulting intersection is empty;
it means that there is a contradiction between the hypotheses of the
logical property and Gappa will use it to prove all the goals of the
logical property.

Once the logical property has been proved, a formal proof is generated
by retracing the paths that were followed when computing the ranges.

\subsection{Hints} \label{sec:hints}
When Gappa is not able to satisfy the goal of the logical property, this
does not necessarily mean that the property is false. It may just mean
that Gappa has not enough information about the expressions it tries to
bound.

It is possible to help Gappa in these situations by providing
\emph{rewriting hints}. These are rewriting rules similar to those
presented above in the case of rounding, but whose usefulness is
specific to the problem at hand.

\subsubsection{Explicit hints}
A hint has the following form:
\begin{lstlisting}[language=gappa,numbers=none,frame=lines]
Expr1 -> Expr2;
\end{lstlisting}

It is used to give the following information to Gappa: ``I believe for some
reason  that, should you need to compute an interval for
\texttt{Expr1}, you might get a tighter interval by trying the mathematically
equivalent \texttt{Expr2}''. This fuzzy formulation is better explained by considering the
following examples.

\begin{enumerate}
\item The ``some reason'' in question will typically be that the
  programmer knows that expressions \texttt{A}, \texttt{B} and
  \texttt{C} are different approximations of the same quantity, and
  furthermore that \texttt{A} is an approximation to \texttt{B} which
  is an approximation to \texttt{C}. As previously, this means that
  these variables are correlated, and the adequate hint to give in
  this case is
\begin{lstlisting}[language=gappa,numbers=none,frame=lines]
A - C -> (A - B) + (B - C);
\end{lstlisting}
It will suggest to Gappa to first compute intervals for \texttt{A-B}
and \texttt{B-C}, and then to sum them to get
an interval for \texttt{A-C}.

As there are an infinite number of arbitrary \texttt{B} expressions that
can be inserted in the right hand side expression, Gappa does not try to
apply every possible rewriting rule when it encounters \texttt{A-C}.
However, as \ref{sec:automatic-hints} will show, Gappa usually infers
some useful \texttt{B} expressions and applies the rewriting rules
automatically.

\item Relative errors can be manipulated similarly. The hint to use in this case is 
\begin{lstlisting}[language=gappa,numbers=none,frame=lines]
(A-C)/C -> (A-B)/B + (B-C)/C + ((A-B)/B)*((B-C)/C);
\end{lstlisting}

This is still a mathematical identity, as one may check easily.
Again, Gappa tries to infer some useful \texttt{B} expressions and to
apply the corresponding rewriting rules.

\item When \texttt{x} is an approximation of $\mathtt{MX}$ and a
  relative error $\epsilon = \frac{\mathtt{x} -
    \mathtt{MX}}{\mathtt{MX}}$ is known by the tool, \texttt{x} can be
  rewritten $\mathtt{MX} \cdot \left( 1 + \epsilon\right)$.  This kind
  of hint is useful in combination with the following one.

\item When manipulating fractional terms such as
  $\frac{\mathtt{Expr1}}{\mathtt{Expr2}}$ where $\mathtt{Expr1}$ and $\mathtt{Expr2}$ are
  correlated (for example one approximating the other), the interval
  division fails to give useful results if the interval for $\mathtt{Expr2}$
  comes close to $0$. In this case, one will try to write
  $\mathtt{Expr1} = \mathtt{A} \cdot \mathtt{Expr3}$ and
  $\mathtt{Expr2} = \mathtt{A} \cdot \mathtt{Expr4}$, so that the
  interval on $\mathtt{Expr4}$ does not come close to $0$ anymore. The
  following hint is then appropriate:
\begin{lstlisting}[language=gappa,numbers=none,frame=lines]
Expr1 / Expr2 -> Expr3 / Expr4;
\end{lstlisting}
This rewriting rule is only valid if $\mathtt{A}$ is not zero, so the
case $\mathtt{A}=0$ has to be handled separately.
So that Gappa does not apply the rule in an invalid context, a constraint
on $\mathtt{A}$ can be added. The rule thus becomes:
\begin{lstlisting}[language=gappa,numbers=none,frame=lines]
Expr1 / Expr2 -> Expr3 / Expr4   { A <> 0 };
\end{lstlisting}
\end{enumerate}

All these hints are correct if both sides are mathematically
equivalent. Gappa therefore checks this automatically. If the test
fails, it emits a warning to the user that he or she
must review the hint by hand. Therefore, writing even
complex hints is very safe: one may not introduce an error in the
proof by writing hints which do not emit warnings. 

Besides, useless hints may consume execution time as Gappa tries in
vain to use them, but if they are useless, they will be silently
ignored in the final proof. Therefore, writing useless hints is
essentially harmless.

\subsubsection{Automatic hints}\label{sec:automatic-hints}

After using Gappa to prove several elementary functions, one thing
became clear: Users kept writing the same hints, typically of the three first kinds
enumerated in \ref{sec:hints}.

Gappa was therefore modified to introduce a new kind of hint:
\begin{lstlisting}[language=gappa,numbers=none,frame=lines]
Expr1 ~ Expr2;
\end{lstlisting}
that reads ``\texttt{Expr1} approximates \texttt{Expr2}''. This has the
effect of introducing rewriting hints both for absolute and relative
differences involving \texttt{Expr1} or \texttt{Expr2}. There may be
useless hints among such automatic hints, but again they will be
mostly harmless.

For instance, when it encounters an expression of the form
\texttt{Expr1 - Expr3} (for any expression \texttt{Expr3}),
Gappa automatically tries the rule \texttt{Expr1-Expr3 -> (Expr1-Expr2) +
(Expr2-Expr3)}. And when it encounters \texttt{Expr3 - Expr2}, it tries
the rewriting rule \texttt{Expr3-Expr2 -> (Expr3-Expr1) + (Expr1-Expr2)}.

By default, Gappa assumes that \texttt{Expr1} approximates \texttt{Expr2}
if \texttt{Expr1} is the rounded value of \texttt{Expr2}. It also makes
this assumption when an enclosure of the absolute or relative error between
\texttt{Expr1} and \texttt{Expr2} appears in the hypotheses of the
logical proposition Gappa has to prove.

Remark: It is very important, in the various differences appearing in
all these expressions, that the least accurate term is written first
and the most accurate written last. The main reason is that the
theorems of Gappa's database apply to expressions written in this order.
This ordering convention prevents a combinatorial explosion on the number
of paths to explore.

\subsubsection{Bisection and dichotomy hints}\label{sec:dicho-hints}
Finally, it is possible to instruct Gappa to split some intervals and perform its exploration on the resulting sub-intervals. 
There are several possibilities. For instance, the following hint
\begin{lstlisting}[language=gappa,numbers=none,frame=lines]
$ z in (-1,2);
\end{lstlisting} 
reads ``Better enclosures may be obtained by separately considering
the sub-cases $\mathtt{z} \le -1$, $-1 \le \mathtt{z} \le 2$, and
$2 \le \mathtt{z}$.''

The following hint finds the splitting points automatically by performing
a dichotomy on the interval of $\mathtt{z}$ until the part of the goal
corresponding to $\mathtt{Expr}$ as been satisfied for all the
sub-intervals of $\mathtt{z}$.
\begin{lstlisting}[language=gappa,numbers=none,frame=lines]
Expr $ z;
\end{lstlisting} 

\subsubsection{Writing hints in an interactive way}

Gappa has evolved to include more and more automatic hints, but most
real-world proofs still require writing complex, problem-specific
hints. Finding the right hint that Gappa needs could be quite complex
and would require completely mastering its theorem database and the
algorithms used by its engine. Fortunately, a much simpler way is to
build the proof incrementally and question the tool by adding and
removing intermediate goals to prove, as the extended example in next
section will show. Before that, we first describe
the outline of the methodology we use to prove elementary functions.

\section{Proving elementary functions using Gappa}
\label{sec:tutorial}

As in every proof work, style is important when working with Gappa: in
a machine-checked proof, bad style will not in principle endanger the
validity of the proof, but it may prevent its author to get to the
end. In the \crlibm{} framework, it may hinder acceptance of
machine-checked proofs among new developers.

Gappa does not impose a style, and when we started using it there was
no previous experience to get inspiration from. After a few months of use, we
had improved our ``coding style'' in Gappa, so that the proofs
were much more concise and readable than the earlier ones. We had also
set up a methodology that works well for elementary functions. This
section is an attempt to describe this methodology and style. We
are aware that they may be inadequate for other applications, and that
even for elementary functions they could be improved further.

The methodology consists in three steps, which correspond to the three sections of a Gappa input file.
\begin{itemize}
\item First, the C code is translated into Gappa equations, in a way
  that ensures that the Gappa proof will indeed prove some property of
  this program (and not of some other similar program). Then
  equations are added describing what the program is supposed to
  implement. Usually, these equations are also in correspondence with
  the code.
\item Then, the property to prove is added.  It is usually in the form \texttt{hypotheses -> properties},
  where the hypotheses are known bounds on the inputs, or contribution
  to the error determined outside Gappa, like the approximation
  errors.
\item Finally, one has to add \emph{hints} to help  Gappa
  complete the proof. This last part is built incrementally.
\end{itemize}
The following sections detail these three steps.

\subsection{Translating a FP program}

We consider again the following C code, where we have added the
constants:
\begin{lstlisting}[language=C,firstnumber=1]
  s3 = -1.6666666666666665741e-01; 
  s5 =  8.3333333333333332177e-03; 
  s7 = -1.9841269841269841253e-04;
  yh2 = yh*yh;					   
  ts = yh2 * (s3 + yh2*(s5 + yh2*s7));	   
  Fast2Sum(sh,sl,   yh,  yl + yh*ts);	           
\end{lstlisting}


There is a lot of rounding operations in this code, so the first thing
to do is to define Gappa rounding operators for the rounding modes
used in the program. In our example, we use the following line to
define \texttt{IEEEdouble} as a shortcut for IEEE-compliant rounding
to the nearest double, which is the mode used in \crlibm.
\begin{lstlisting}[language=gappa,numbers=none,frame=lines]
@IEEEdouble = float<ieee_64,ne>;
\end{lstlisting}

Then, if the C code is itself sufficiently simple and clean, the
translation step only consists in making explicit the rounding
operations implicit in the C source code.  To start with, the
constants \texttt{s3}, \texttt{s5} and \texttt{s7} are given as
decimal strings, and the C compilers we use convert them to (binary)
double-precision FP numbers with round to nearest. We ensure that
Gappa works with these same constants as the compiled C code by
inserting explicit rounding operations:
\begin{lstlisting}[language=gappa,numbers=none,frame=lines]
  s3 = IEEEdouble(-1.6666666666666665741e-01); 
  s5 = IEEEdouble( 8.3333333333333332177e-03); 
  s7 = IEEEdouble(-1.9841269841269841253e-04);
\end{lstlisting}

Then we have to do the same for all the roundings hidden behind C arithmetic operations.   Adding by hand all the rounding operators, however, would be tedious
and error-prone, and would make the Gappa syntax so different from the
C syntax that it would degrade confidence and maintainability.
Besides, one would have to apply without error the rules (well
specified by the C99 standard \cite{C99}) governing for instance implicit parentheses in a C
expression. For these reasons, Gappa has a syntax that instructs it to
perform this task automatically. The following Gappa lines
\begin{lstlisting}[language=gappa,numbers=none,frame=lines]
  yh2 IEEEdouble=  yh * yh;
  ts  IEEEdouble=  yh2 * (s3 + yh2*(s5 + yh2*s7));
\end{lstlisting}
define the same mathematical relation between their right-hand side
and left-hand side as the corresponding lines of the C programs. This,
of course, is only true under the following conditions:
\begin{itemize}
\item all the C variables are double-precision variables,
\item the Gappa variables on the right-hand side represent them,
\item the compiler/OS/processor combination used to
  process the C code respects the C99 and IEEE-754 standards and
  computes in double-precision arithmetic.
\end{itemize}

Finally, we have to express in Gappa the Fast2Sum algorithm. Where in
C it is a macro or function call, for our purpose we prefer to ignore
this complexity and simply express in Gappa the resulting behaviour,
which is a sum without error (we have here to trust an external proof
of this behaviour \cite{Dekker71,Knuth97}):

\begin{lstlisting}[language=gappa,numbers=none,frame=lines]
  r   IEEEdouble=  yl + yh*ts;
  s             =  yh + r;   #  the Fast2Sum is exact. s is sh+sl 
\end{lstlisting}

  Note that we are interested in the relative error of the sum
\texttt{s} with respect to the exact sine, and for this purpose the
fact that \texttt{s} has to be represented as a sum of two doubles in
C is irrelevant.

More importantly, this adds another condition for this code
translation to be faithful: As the proof of the Fast2Sum has as hypothesis that the exponent
of \texttt{yh} is larger than or equal to that of \texttt{r}, we now
have to prove that. We will simply add this goal to the theorem  to prove.


As a summary,  for straight-line program segments with mostly
double-precision variables, a set of corresponding Gappa definitions
can be obtained straightforwardly by just replacing the C \texttt{=}
with Gappa \texttt{IEEEdouble=}, a very safe operation. 


\subsection{Defining ideal values} \label{subsec:ideal}

To analyse this code, we now need a ``mathematically ideal'' definition of
all the variables, a reference with respect to which the error is
computed.  This notion of mathematically ideal may be quite
subtle: What is the mathematically ideal of \texttt{yh2}? It could be
\begin{itemize}
\item  the exact square of \texttt{yh} (without rounding), or
\item  the exact square of
  \texttt{yh+yl} that \texttt{yh} approximates, or
\item  the exact square of
  the ideal reduced argument $y$ , which is usually irrational.
\end{itemize}

The right choice depends on the properties to prove. Here, the mathematically ideal value for both \texttt{yh+yl} and
\texttt{yh} will be the ideal reduced argument, which we note
\texttt{My}. Similarly, the purest mathematical value that
\texttt{yh2} approximates is noted \texttt{My2} and will be defined as
\texttt{My2 = My*My}. 

For the polynomial approximation
\texttt{ts}, we could choose, as mathematical ideal, either the value of the
function that the polynomial approximates, or the value of the same
polynomial, but computed on \texttt{My} and without rounding error.
Here we chose the latter, although it is the lesser ideal.

Here come a few naming conventions. The first was obviously that the
Gappa variables that mimick C variables have the same name. We also impose 
the convention that such variables begin with a lowercase letter. In addition, Gappa variables for
mathematically ideal terms will begin with a ``\texttt{M}''. The
other intermediate Gappa variables should begin with capital letters
to distinguish them from variables mimicking the code. Of course,
related variables should have related and, wherever possible, explicit
names.  Again, these are conventions and are part of a proof style,
not part of Gappa syntax: the capitalization will give no information
to the tool, and neither will the fact that variables have related
names.

For instance, it will be convenient to define a variable equal to \texttt{yh+yl}:
\begin{lstlisting}[language=gappa,numbers=none,frame=lines]
Yhl = yh + yl;
\end{lstlisting}

\subsection{Defining what the code is supposed to compute}

Defining mathematically ideal values resumes to  defining in Gappa what the C code is
supposed to implement. For instance, using our previous conventions, the line for \texttt{ts} was
probably evaluating the value of the same polynomial of the ideal
\texttt{My}:
\begin{lstlisting}[language=gappa,numbers=none,frame=lines]
My2 = My*My;
Mts = My2 * (s3 + My2*(s5 + My2*s7));
\end{lstlisting}

We have kept the polynomial coefficients in lower case: As already
discussed in Section~\ref{sec:bound}, the polynomial thus defined
nevertheless belongs to the set of polynomial with real coefficients,
and we have means  to compute (outside Gappa) a bound of its relative error with
respect to the function it approximates.

The link between \texttt{ts} and the polynomial approximating the sine
is also best expressed using  mathematically ideal values:

\begin{lstlisting}[language=gappa,numbers=none,frame=lines]
PolySinY = My + My*Mts;
\end{lstlisting}

To sum up, \texttt{PolySinY} is the actual polynomial with the same
coefficients \texttt{s3} to \texttt{s7} as in the C code, but
evaluated without rounding error, and evaluated on the ideal value
that \texttt{yh+yl} approximates.

Another approach could be to use a Taylor polynomial,
in which case the approximation error would be given by the rest in
the Taylor formula, the ideal polynomial would be the Taylor one, it
would have ideal Taylor coefficients (beginning with M), some of which
would have to be rounded to FP numbers to appear in the program
(lowercase). Gappa could handle it, too, but it would be less convenient.

Another crucial question is, how do we define the real, ideal, mathematical
function which we eventually approximate? Gappa has no builtin sine or
logarithm. The current approach can be described in English as:
``$\sin(\mathtt{My})$ is a value which, as long as $\mathtt{My}$ is smaller than \texttt{6.29e-3}, remains within
a relative distance of \texttt{2.26e-24} of our ideal polynomial''. In Gappa,
this will translate to some hypotheses in the property to prove:
\begin{lstlisting}[language=gappa,numbers=none,frame=lines]
    |Yhl| in [0, 6.29e-03]
 /\ |(PolySinY - SinY)/SinY| <= 2.26e-24
 /\ ... # (more hypotheses, see below)
->    epstotal in ?
\end{lstlisting}

Here the interval of \texttt{Yhl} is defined by the range reduction,
and its bound has been computed separately: it is $\pi/512$, plus some
margin that accounts for the inexactness of argument reduction.

Similarly, the relative distance between the sine and the polynomial
on this interval is computed outside Gappa. We used to use Maple's
infinite norm, but it only returns an approximation,
so this was in principle a weakness of the proof. We now use a safer,
interval-based approach \cite{ChevLaut2007QSIC}, implemented in the
Sollya tool\footnote{\url{http://sollya.gforge.inria.fr/}}.

Of course this tool provides a theorem which could be expressed as\\
\verb!|Yhl| in [0, 6.29e-03] -> |(PolySinY - SinY)/SinY| <= 2.26e-24!. We just inject the property of this theorem as an hypothesis in Gappa.

Concerning style, it will be more convenient to have a variable
defined as this approximation error. Similarly, it will make the proof
much clearer and concise to add, from the beginning, as many
definitions as possible for the various terms and errors involved in
the computation:
\begin{lstlisting}[language=gappa,numbers=none,frame=lines]
Epsargred = (Yhl - My)/My;           # argument reduction error
Epsapprox = (PolySinY - SinY)/SinY;  # polynomial approximation error
Epsround = (s - PolySinY)/PolySinY;  # rounding errors in the polynomial evaluation
Epstotal = (s - SinY)/SinY;          # total error
\end{lstlisting}

\subsection{Defining the property to prove}

With the previous introduction, the theorem to prove, expressed as
implications using classical first-order logic, is stated as follows:

\begin{lstlisting}[language=gappa,numbers=none,frame=lines]
{
  # Hypotheses 
    |Yhl| in [0, 6.29e-03] 
 /\ |Epsargred| <= 2.53e-23
 /\ |Epsapprox| <= 2.26e-24

-> 

   Epstotal in ?     # the main goal of our theorem
/\ |r/yh| in [0,1]   # the condition for the  Fast2Sum to be valid 
}
\end{lstlisting}

The full initial Gappa script is given below. It adds a more accurate
definition of \texttt{yh} and \texttt{yl}, stating that they are
double-precision numbers and that they form a disjoint double-double.
It also adds a lower bound on the absolute value of the reduced
argument, obtained thanks to Kahan/Douglas algorithm. This lower bound
is important because it will enable most values to stay away from
zero, which ensures that relative errors are not arbitrarily big due to
underflow.

\begin{lstlisting}[language=gappa,frame=lines,firstnumber=1,label=lst:gappainit,caption=The initial Gappa file.]
@IEEEdouble = float<ieee_64,ne>;

# yh+yl is a double-double (call it Yhl)
yh = IEEEdouble(dummy1);
yl = IEEEdouble(dummy2);
Yhl = yh + yl;       # Below, there is also an hypothesis stating that yl<ulp(yh)

#--------------- Transcription of the C code --------------------------

s3 = IEEEdouble(-1.6666666666666665741e-01); 
s5 = IEEEdouble( 8.3333333333333332177e-03); 
s7 = IEEEdouble(-1.9841269841269841253e-04);
yh2 IEEEdouble=  yh * yh;
ts  IEEEdouble=  yh2 * (s3 + yh2*(s5 + yh2*s7));
r   IEEEdouble=  yl + yh*ts;
s             =  yh + r;   # no rounding, it is the Fast2Sum 

#-------- Mathematical definition of what we are approximating --------

My2 = My*My;
Mts = My2 * (s3 + My2*(s5 + My2*s7));
PolySinY = My + My*Mts;

Epsargred = (Yhl - My)/My;           # argument reduction error
Epsapprox = (PolySinY - SinY)/SinY;  # polynomial approximation error
Epsround = (s - PolySinY)/PolySinY;  # rounding errors in the polynomial evaluation
Epstotal = (s - SinY)/SinY;          # total error

#----------------------  The theorem to prove --------------------------
{
  # Hypotheses 
    |yl / yh| <= 1b-53
 /\ |Yhl| in [1b-200, 6.29e-03] # lower bound guaranteed by Kahan-Douglas algorithm
 /\ |Epsargred| <= 2.53e-23
 /\ |Epsapprox| <= 2.26e-24

->

#goal to prove
   Epstotal in ?     # [-1b-67, 1b-67]
/\ |r/yh| <= 1
}
\end{lstlisting}

\subsection{With a little help from the user}

Invoking Gappa on this file produces the following output:
\begin{small}
\begin{verbatim}
Warning: no path was found for Epstotal.
Warning: no path was found for |r / yh|.

Results for |yl / yh| in [0, 1.11022e-16] and |Yhl| in [6.22302e-61, 0.00629] 
        and |Epsargred| in [0, 2.53e-23] and |Epsapprox| in [0, 2.26e-24]:
Warning: some enclosures were not satisfied.
\end{verbatim}
\end{small}
This means that Gappa needs some help, in the form of hints. Where to
start? There are several way to interact with the tool to understand
where it fails.

\begin{itemize}

\item We may add additional goals to obtain enclosures for intermediate variables. For instance, adding the goal \verb!|My| in ?!, we obtain the following answer 
  \begin{small}
\begin{verbatim}
|My| in [0, 0.00629]
\end{verbatim}
  \end{small}
  Gappa was able to deduce this enclosure from the enclosure of
  \texttt{Yhl} (hypothesis) and the definition of \texttt{Epsargred}.
  Similarly, we may check for instance that the built-in engine is able to build a
  good enclosure of \texttt{PolySinY}, but not of \texttt{s}.

\item We may add additional hypotheses and see what progress they
  entail. For instance, providing a dummy \texttt{Epsround} as an
  hypothesis allows Gappa to complete the proof, thanks to its
  automatic hints.
\end{itemize}

This way it is possible to track the point where Gappa's engine gets
lost, and provide hints to help it.


In our case, the best thing to do is to express all the approximation
layers detailed in Section~\ref{sec:troislignes}.
Written as Gappa equations, we get:

 \begin{lstlisting}[numbers=none,frame=lines]
# Layers of approximation on s
S1 = yh + (yl + IEEEdouble(yh*ts));   # s without last rounding
S2 = yh + (yl + yh*ts);               # removing penultimate rounding, too
S3 = (yh+yl) + (yh+yl)*ts;            # putting back yl which was neglected

Eps1 = (s-S1)/S1;
Eps2 = (S1-S2)/S2;
Eps3 = (S2-S3)/S3;
Eps4 = (S3-PolySinY)/PolySinY;
\end{lstlisting}

Remark again that all these relative errors are defined relatively to the most accurate term. 

We may add goals for these new relative errors: Gappa will be unable to bound any of them. We have to provide hints. 

Consider only \texttt{Eps4}, the relative difference between
\texttt{S3} and \texttt{PolySinY} -- which we saw Gappa is able to
bound.  Both \texttt{S3} and \texttt{PolySinY} are polynomial
expressions without any rounding, and with identical coefficients.
Therefore, the difference between them resumes to the difference
between \texttt{Yhl=yh+yl}, used in \texttt{S3}, and \texttt{My}, used
in \texttt{PolySinY}. We precisely have a measure of this difference:
it is \texttt{Epsargred}. The hint we have to provide to Gappa should
therefore express \texttt{Eps4} as a function of \texttt{Epsargred}
which, when evaluated by intervals, will provide a tight enclosure.
Here is a generic technique to obtain such an hint.  We start with
\texttt{Eps4 -> (S3-PolySinY)/PolySinY}, which is just the definition
of \texttt{Eps4}, and we rewrite it incrementally until we have
obtained an expression involving \texttt{Epsargred}.  In the following
hint, we have left, for the purpose of this tutorial, the intermediate
rewriting steps commented out.

 \begin{lstlisting}[numbers=none,frame=lines]
Eps4 -> 
#  (S3-PolySinY)/PolySinY;
#   S3/PolySinY - 1;
#   ((yh+yl) + (yh+yl)*ts) / (My + My*Mts)  - 1;
#   ((yh+yl)/My) * (1+ts)/(1+Mts)  - 1;
#   (Epsargred+1) * (1+ts)/(1+Mts)  - 1;
#    Epsargred * (1+ts)/(1+Mts)   +  1 * (1+ts)/(1+Mts)  - 1;
#    Epsargred * (1+ts)/(1+Mts)   +  (ts-Mts)/(1+Mts);
   Epsargred * (1+ts)/(1+Mts)   +  Mts*((ts-Mts)/Mts) / (1+Mts);
\end{lstlisting}

Considering the orders of magnitudes (see
Section~\ref{sec:troislignes}), a naive interval evaluation of this
last expression will be very accurate. Indeed, \texttt{Mts} as well as
\texttt{ts} are very small compared to 1, therefore the first term is
close to \texttt{Epsargred}. The second is the relative error of
\texttt{ts} with respect to \texttt{Mts} (expected to be no larger
than $2^{-52}$), multiplied by \texttt{Mts} which is smaller than
$2^{-14}$. We therefore have a sum of two small terms which should
provide a small enclosure.

Still, even with this hint, Gappa still fails to provide an enclosure
for \texttt{Eps4}. Adding the goals \texttt{Mts in ?} and \texttt{ts
  in ?}, we observe that \texttt{Mts} is properly enclosed, but not
\texttt{ts}. We therefore add a definition of the relative error of
\texttt{ts} with respect to \texttt{Mts}:

\begin{lstlisting}[numbers=none,frame=lines]
 EpstsMts = (ts-Mts)/Mts ;
\end{lstlisting}
and we perform the same analysis: we describe the succession of
approximation layers between \texttt{ts} and \texttt{Mts}, define
intermediate error terms for them, and provide hints for bounding
them. This will in turn require to explain to Gappa how to go from
\texttt{My} to \texttt{yh2} through three approximation layers.

We will not detail this process line by line. The final Gappa script 
is given in appendix, and 
is available from the distribution of \crlibm
\footnote{http://lipforge.ens-lyon.fr/www/crlibm/}.

\subsection{Summing up}

Writing hints is the most time-consuming part of the proof, because it
is the part where the designer's intelligence is required. However, we
hope to have shown that it may be done very incrementally.

The example chosen in this article is actually quite complex: its Gappa
proof consists of more than 150 lines, half of which are hints. The
bound found on \texttt{Epstotal} is $2^{-67,24}$ and is obtained in a
few seconds on a recent machine (the time can be longer when there is
a dichotomy).  The resulting Coq proof is more than 7000 lines long.

Some functions are simpler. We could write the proof of a logarithm
implementation \cite{DinLauMul2007:log} with a few hints only~%
\cite{DinLauMelq2005LIP:gappa}. One reason is that the logarithm never
comes close to 0, so the full proof can be handled only with absolute
errors, for which writing hints is much lighter.


\section{Conclusion and perspectives}




Validating tight error bounds on the low-level, optimized
floating-point code typical of elementary functions has always been a
challenge, as many sources of errors cumulate their effect.
Gappa is a high-level proof assistant that is well suited to this kind of
proofs. 

Using Gappa, it is easy to translate a part of a C program into a
mathematical description of the operations involved with fair
confidence that this translation is faithful. Expressing implicit
mathematical knowledge one may have about the code and its context is
also easy. Gappa uses interval arithmetic to manage the ranges and
errors involved in numerical code. It handles most of the
decorrelation problems automatically thanks to its built-in rewriting
rules, and an engine which explores the possible rewriting of
expressions to minimize the size of the intervals. If decorrelation
remains, Gappa allows one to provide new rewriting rules, but checks
them. All this is well founded on a library of theorems which allow the
obtained computation to be translated to a proof checkable by a
lower-level proof assistant such as Coq and PVS. 
Finally, the tool can be questioned during the process of building the
proof so that this process may be conducted interactively. 

Therefore, it is possible to get quickly a fully validated proof with
good confidence that this proof indeed proves property of the initial
code. Gappa is by no means automatic: to apply it on a given piece of
code requires exactly the same {knowledge} and cleverness a paper
proof would. However, it requires much less {work}.

The current \crlibm{} distribution contains several bits of
proofs using Gappa at several stages of its development. Although this
development is not over, the current version (0.9) is very stable
and we may safely consider generalizing the use of this tool in the
future developments of \crlibm{}. It also took 6 months to
develop a methodology and style well suited to the validation of
elementary functions. This paper presented this aspect as well. Very
probably, new problems will arise as we try to apply this methodology
to new functions, so that it will need to be refined further. 

Iterative codes are currently out of scope of our
  methodology, although it could be used for instance to
  prove loop invariants.


\appendix

\section{Complete Gappa script}

\begin{lstlisting}[language=gappa,frame=lines,firstnumber=1,label=lst:gappainit,caption=The complete Gappa file.]
# test with   gappa -Munconstrained < sin.gappa
# The proof is not complete, as it doesn't work without -Munconstrained. 
# What it means is that Gappa is unable to prove that some denominators are not null.
# It's OK for practical purposes, but it takes some more work to get a formal proof. 

@IEEEdouble = float<ieee_64,ne>;
# Convention 1: uncapitalized variables match the variables in the C code. Other variables begin with a capital letter
# Convention 2: variables beginning with "M" are mathematical ideal

# yh+yl is a double-double (call it Yhl)

yh = IEEEdouble(dummy1);
yl = IEEEdouble(dummy2);
Yhl = yh + yl; # There is also an hypothesis stating that yl<ulp(yh)

#--------------- Transcription of the C code --------------------------

s3 = IEEEdouble(-1.6666666666666665741e-01); 
s5 = IEEEdouble( 8.3333333333333332177e-03); 
s7 = IEEEdouble(-1.9841269841269841253e-04);

yh2 IEEEdouble=  yh * yh;
ts  IEEEdouble=  yh2 * (s3 + yh2*(s5 + yh2*s7));
r   IEEEdouble=  yl + yh*ts;
s             =  yh + r;   # no rounding, it is the Fast2Sum 

#-------- Mathematical definition of what we are approximating --------

My2 = My*My;
Mts = My2 * (s3 + My2*(s5 + My2*s7));
PolySinY = My + My*Mts;

Epsargred = (Yhl - My)/My;           # argument reduction error
Epsapprox = (PolySinY - SinY)/SinY;  # polynomial approximation error
Epsround = (s - PolySinY)/PolySinY;  # rounding errors in the polynomial evaluation
Epstotal = (s - SinY)/SinY;          # total error



# Layers of approximation on s
S1 = yh + (yl + IEEEdouble(yh*ts));   # remove last round
S2 = yh + (yl + yh*ts);               # remove penultimate round
S3 = (yh+yl) + (yh+yl)*ts;            # put yl back in

Eps1 = (s-S1)/S1;
Eps2 = (S1-S2)/S2;
Eps3 = (S2-S3)/S3;
Eps4 = (S3-PolySinY)/PolySinY;


yhts = IEEEdouble(yh*ts);             # just to make the hints lighter
p3 IEEEdouble=  s3 + yh2*(s5 + yh2*s7);  # idem


tsNoRound = yh2 * (s3 + yh2*(s5 + yh2*s7));

# A few definitions mostly to benefit from automatic hints.
EpstsMts = (ts-Mts)/Mts;
EpstsNoRoundMts = (tsNoRound - Mts)/Mts;

Epsy2 = (yh2-My2)/My2;
Epsy2_argred = (Yhl*Yhl-My2)/My2;
Epsy2_negl_yl = (yh*yh-Yhl*Yhl)/(Yhl*Yhl);
Epsy2_rnd = (yh2-yh*yh)/(yh*yh);


#----------------------  The theorem to prove --------------------------
{
  # Hypotheses 
    |yl / yh| <= 1b-53
 /\ |Yhl| in [1b-200, 6.29e-03] # lower bound guaranteed by Kahan-Douglas algorithm
 /\ |yh| in [1b-1000, 1] # some huge range for ensuring that yh is not zero
 /\ |Epsargred| <= 2.53e-23
 /\ |Epsapprox| <= 2.26e-24

->

#goal to prove
   Epstotal in ?# [-1b-67, 1b-67]
/\ |r/yh| <= 1
#/\ |My| in [1b-400, 6.29e-03]
}

# ---------------------- Hints ----------------------------------
$ Yhl in (0);

# First, the hints for Epsround

s~S1;
S1~S2;
S2~S3;
S3~PolySinY;

Eps4 -> #  (S3-PolySinY)/PolySinY;
#   S3/PolySinY - 1;
#   ((yh+yl) + (yh+yl)*ts) / (My + My*Mts)  - 1;
#   ((yh+yl)/My) * (1+ts)/(1+Mts)  - 1;
#   (Epsargred+1) * (1+ts)/(1+Mts)  - 1;
#    Epsargred * (1+ts)/(1+Mts)   +  1 * (1+ts)/(1+Mts)  - 1;
#    Epsargred * (1+ts)/(1+Mts)   +  (ts-Mts)/(1+Mts);
   Epsargred * (1+ts)/(1+Mts)   +  Mts*((ts-Mts)/Mts) / (1+Mts);

# Now we just need to bound ts-Mts: 
ts ~ tsNoRound;
(tsNoRound - Mts)/Mts -> 
#     yh2/My2  * (s3 + yh2*(s5 + yh2*s7)) / (s3 + My2*(s5 + My2*s7))  - 1 ;
      (1+Epsy2)  * (s3 + yh2*(s5 + yh2*s7)) / (s3 + My2*(s5 + My2*s7)) -1;
# Now we just need to express My2 in terms of yh2 and Epsy2
My2 -> yh2/(1+Epsy2);

yh ~ Yhl;
(yh - Yhl) / Yhl -> 1 / (1 + yl / yh) - 1;

Eps3 -> 
# (S2-S3)/S3 
#  S2/S3 - 1; 
#   (yh + (yl + yh*ts)) / ((yh+yl) + (yh+yl)*ts)   - 1 ;
#   ((yh+yl) + (yh+yl)*ts - yl*ts) / ((yh+yl) + (yh+yl)*ts)   - 1 ;
#    - yl*ts / ((yh+yl) + (yh+yl)*ts)   ;
#    - (yl/Yhl)  * (ts / (1+ts))   ;
     ((yh-Yhl)/Yhl)  * (ts / (1+ts))   ; # change sign to have the expression of a rounding error


Eps2 -> # (S1-S2)/S2; 
#   (yh + (yl + IEEEdouble(yh*ts))) / (yh + (yl + yh*ts))   -1 ;
#   (IEEEdouble(yh*ts) - yh*ts) / (yh + yl + yh*ts) ;
#   ((IEEEdouble(yh*ts) - yh*ts)/(yh*ts)) / ( (yh+yl)/(yh*ts) + 1 ) ;
    ts * ((IEEEdouble(yh*ts) - yh*ts)/(yh*ts)) / ( 1 + yl/yh + ts ) ;

yhts/yh -> ts*((yhts-yh*ts)/(yh*ts) + 1);


(yl+yhts)/yh -> yl/yh + yhts/yh;
\end{lstlisting}

\bibliographystyle{IEEEtran.bst} 
\bibliography{arith}

\begin{thebibliography}{10}
\providecommand{\url}[1]{#1}
\csname url@samestyle\endcsname
\providecommand{\newblock}{\relax}
\providecommand{\bibinfo}[2]{#2}
\providecommand{\BIBentrySTDinterwordspacing}{\spaceskip=0pt\relax}
\providecommand{\BIBentryALTinterwordstretchfactor}{4}
\providecommand{\BIBentryALTinterwordspacing}{\spaceskip=\fontdimen2\font plus
\BIBentryALTinterwordstretchfactor\fontdimen3\font minus
  \fontdimen4\font\relax}
\providecommand{\BIBforeignlanguage}[2]{{%
\expandafter\ifx\csname l@#1\endcsname\relax
\typeout{** WARNING: IEEEtran.bst: No hyphenation pattern has been}%
\typeout{** loaded for the language `#1'. Using the pattern for}%
\typeout{** the default language instead.}%
\else
\language=\csname l@#1\endcsname
\fi
#2}}
\providecommand{\BIBdecl}{\relax}
\BIBdecl

\bibitem{IEEE754}
ANSI/IEEE, \emph{Standard 754-1985 for Binary Floating-Point Arithmetic (also
  {IEC} 60559)}, 1985.

\bibitem{Goldberg91}
D.~Goldberg, ``What every computer scientist should know about floating-point
  arithmetic,'' \emph{ACM Computing Surveys}, vol.~23, no.~1, pp. 5--47, Mar.
  1991.

\bibitem{Moore66}
R.~E. Moore, \emph{Interval analysis}.\hskip 1em plus 0.5em minus 0.4em\relax
  Prentice Hall, 1966.

\bibitem{DinErshGast2005}
\BIBentryALTinterwordspacing
F.~de~Dinechin, A.~Ershov, and N.~Gast, ``Towards the post-ultimate libm,'' in
  \emph{17th Symposium on Computer Arithmetic}.\hskip 1em plus 0.5em minus
  0.4em\relax IEEE Computer Society Press, Jun. 2005, pp. 288--295. [Online].
  Available:
  \url{http://perso.ens-lyon.fr/florent.de.dinechin/recherche/publis/2005-Arit%
h.pdf}
\BIBentrySTDinterwordspacing

\bibitem{DinLauMul2007:log}
F.~de~Dinechin, C.~Q. Lauter, and J.-M. Muller, ``Fast and correctly rounded
  logarithms in double-precision,'' \emph{Theoretical Informatics and
  Applications}, vol.~41, pp. 85--102, 2007.

\bibitem{Har97}
J.~Harrison, ``Floating point verification in {HOL} light: the exponential
  function,'' University of Cambridge Computer Laboratory, Technical Report
  428, 1997.

\bibitem{harrison-fmcad2000}
------, ``Formal verification of floating point trigonometric functions,'' in
  \emph{Formal Methods in Computer-Aided Design: Third International Conference
  {FMCAD} 2000}, ser. Lecture Notes in Computer Science, vol. 1954.\hskip 1em
  plus 0.5em minus 0.4em\relax Springer-Verlag, 2000, pp. 217--233.

\bibitem{Dekker71}
T.~J. Dekker, ``A floating point technique for extending the available
  precision,'' \emph{Numerische Mathematik}, vol.~18, no.~3, pp. 224--242,
  1971.

\bibitem{Knuth97}
D.~Knuth, \emph{The Art of Computer Programming, vol.2: Seminumerical
  Algorithms}, 3rd~ed.\hskip 1em plus 0.5em minus 0.4em\relax Addison Wesley,
  1997.

\bibitem{DauRidThe01}
M.~Daumas, L.~Rideau, and L.~Th{\'e}ry, ``A generic library of floating-point
  numbers and its application to exact computing,'' in \emph{Theorem Proving in
  Higher Order Logics: 14th International Conference, TPHOLs 2001}, ser. LNCS,
  R.~J. Boulton and P.~B. Jackson, Eds., vol. 2152.\hskip 1em plus 0.5em minus
  0.4em\relax Springer-Verlag, 2001, pp. 169--184.

\bibitem{Gal86}
S.~Gal, ``Computing elementary functions: A new approach for achieving high
  accuracy and good performance,'' in \emph{Accurate Scientific Computations,
  LNCS 235}.\hskip 1em plus 0.5em minus 0.4em\relax Springer Verlag, 1986, pp.
  1--16.

\bibitem{Tang91}
P.~T.~P. Tang, ``Table lookup algorithms for elementary functions and their
  error analysis,'' in \emph{10th IEEE Symposium on Computer Arithmetic}.\hskip
  1em plus 0.5em minus 0.4em\relax IEEE, Jun. 1991.

\bibitem{StoTan99}
S.~Story and P.~Tang, ``New algorithms for improved transcendental functions on
  {IA}-64,'' in \emph{14th {IEEE} Symposium on Computer Arithmetic}.\hskip 1em
  plus 0.5em minus 0.4em\relax IEEE, Apr. 1999, pp. 4---11.

\bibitem{Markstein2001}
R.-C. Li, P.~Markstein, J.~P. Okada, and J.~W. Thomas, ``The libm library and
  floating-point arithmetic for {HP-UX} on {I}tanium,'' {H}ewlett-{P}ackard
  company, Tech. Rep., april 2001.

\bibitem{Muller2006}
J.-M. Muller, \emph{Elementary Functions, Algorithms and Implementation},
  2nd~ed.\hskip 1em plus 0.5em minus 0.4em\relax Birkh\"{a}user, 2006.

\bibitem{Priest97}
D.~Priest, ``Fast table-driven algorithms for interval elementary functions,''
  in \emph{13th IEEE Symposium on Computer Arithmetic}.\hskip 1em plus 0.5em
  minus 0.4em\relax IEEE, 1997, pp. 168--174.

\bibitem{Markstein2000}
P.~Markstein, \emph{{IA-64} and Elementary Functions: Speed and Precision},
  ser. Hewlett-Packard Professional Books.\hskip 1em plus 0.5em minus
  0.4em\relax Prentice Hall, 2000, iSBN: 0130183482.

\bibitem{ChevLaut2007QSIC}
S.~Chevillard and C.~Lauter, ``A certified infinite norm for the implementation
  of elementary functions,'' in \emph{Seventh International Conference on
  Quality Software ({QSIC} 2007)}.\hskip 1em plus 0.5em minus 0.4em\relax
  {IEEE}, 2007, pp. 153--160.

\bibitem{Ste74}
P.~H. Sterbenz, \emph{Floating point computation}.\hskip 1em plus 0.5em minus
  0.4em\relax Englewood Cliffs, NJ: Prentice-Hall, 1974.

\bibitem{DinMaid2006:RNC}
F.~de~Dinechin and S.~Maidanov, ``Software techniques for perfect elementary
  functions in floating-point interval arithmetic,'' in \emph{Real Numbers and
  Computers}, Jul. 2006.

\bibitem{crlibmweb}
``{CR-Libm}, a library of correctly rounded elementary functions in
  double-precision,'' \url{http://lipforge.ens-lyon.fr/www/crlibm/}.

\bibitem{harrison2003-sqrt}
J.~Harrison, ``Formal verification of square root algorithms,'' \emph{Formal
  Methods in Systems Design}, vol.~22, pp. 143--153, 2003.

\bibitem{filib98}
W.~Hofschuster and W.~Kr\"amer, ``{FI\_LIB}, eine schnelle und portable
  {Funktionsbibliothek} f\"ur reelle {Argumente} und reelle {Intervalle} im
  {IEEE}-double-{Format},'' {Institut f\"ur Wissenschaftliches Rechnen und
  Mathematische Modellbildung, Universit\"at Karlsruhe}, Tech. Rep. Nr. 98/7,
  1998.

\bibitem{filib2005}
W.~Hofschuster, W.~Krämer, M.~Lerch, G.~Tischler, and J.~Wolff~v. Gudenberg,
  ``{filib++} a fast interval library supporting containment computations,''
  \emph{Transactions on Mathematical Software}, 2005, to appear.

\bibitem{DauMel04}
M.~Daumas and G.~Melquiond, ``Generating formally certified bounds on values
  and round-off errors,'' in \emph{6th Conference on Real Numbers and
  Computers}, 2004.

\bibitem{MelPio05}
G.~Melquiond and S.~Pion, ``Formal certification of arithmetic filters for
  geometric predicates,'' in \emph{Proceedings of the 15th IMACS World Congress
  on Computational and Applied Mathematics}, 2005.

\bibitem{Bou97}
S.~Boutin, ``Using reflection to build efficient and certified decision
  procedures,'' in \emph{Third International Symposium on Theoretical Aspects
  of Computer Software}, 1997, pp. 515--529.

\bibitem{HueKahPau04}
\BIBentryALTinterwordspacing
G.~Huet, G.~Kahn, and C.~Paulin-Mohring, \emph{The {Coq} proof assistant: a
  tutorial: version 8.0}, 2004. [Online]. Available:
  \url{ftp://ftp.inria.fr/INRIA/coq/current/doc/Tutorial.pdf.gz}
\BIBentrySTDinterwordspacing

\bibitem{Har00}
\BIBentryALTinterwordspacing
J.~Harrison, \emph{The {HOL} {L}ight manual}, 2000, version~1.1. [Online].
  Available: \url{http://www.cl.cam.ac.uk/users/jrh/hol-light/manual-1.1.pdf}
\BIBentrySTDinterwordspacing

\bibitem{DauMelMun05}
\BIBentryALTinterwordspacing
M.~Daumas, G.~Melquiond, and C.~Mu{\~n}oz, ``Guaranteed proofs using interval
  arithmetic,'' in \emph{Proceedings of the 17th IEEE Symposium on Computer
  Arithmetic}, Cape Cod, Massachusetts, USA, 2005, pp. 188--195. [Online].
  Available:
  \url{http://perso.ens-lyon.fr/guillaume.melquiond/doc/05-arith17-article.pdf}
\BIBentrySTDinterwordspacing

\bibitem{OwrRusSha92}
S.~Owre, J.~M. Rushby, , and N.~Shankar, ``{PVS}: a prototype verification
  system,'' in \emph{11th International Conference on Automated
  Deduction}.\hskip 1em plus 0.5em minus 0.4em\relax Springer, 1992, pp.
  748--752.

\bibitem{C99}
ISO/IEC, \emph{International Standard {ISO/IEC} 9899:1999(E). Programming
  languages -- C}, 1999.

\bibitem{DinLauMelq2005LIP:gappa}
F.~de~Dinechin, C.~Q. Lauter, and G.~Melquiond, ``Assisted verification of
  elementary functions,'' LIP, Tech. Rep. RR2005-43, Sep. 2005.

\end{thebibliography}

\end{document}